\documentclass[11pt,a4paper]{article}

\usepackage[T1]{fontenc}
\usepackage{amsmath,amssymb,amsthm}
\usepackage{geometry}
\usepackage{booktabs}
\usepackage{array}
\usepackage{graphicx}
\usepackage{hyperref}
\newtheorem*{acknowledgement}{Acknowledgements}

\usepackage{amsmath,amssymb,amsthm}

\DeclareMathOperator{\Sym}{Sym}
\DeclareMathOperator{\Aut}{Aut}
\DeclareMathOperator{\Orb}{Orb}

\title{Coherent Configuration on $3$-uniform Hypergraph over Strongly Regular Graph}
\author{
	FAJRI MAULANA
}
\date{}

\newtheorem{theorem}{Theorem}
\newtheorem{proposition}{Proposition}
\newtheorem{corollary}{Corollary}
\newtheorem{definition}{Definition}
\usepackage{enumitem}
\usepackage{longtable}
\usepackage{booktabs}

\begin{document}
	
	\maketitle
	
	\begin{abstract}
This paper investigates coherent configurations on 3-uniform hypergraphs constructed from strongly regular graphs. We establish that the hypergraph preserves the automorphism group of the original graph and induces a nonhomogeneous Schurian coherent configuration. Using the case of Petersen graph, we obtain a configuration with two fibers and thirteen basic relations. We further determine their valencies and nonzero intersection numbers, together with the corresponding adjacency matrix multiplication identities. \\
$Keywords$: Strongly regular graph; 3-uniform hypergraph; coherent configuration; Petersen graph; 
\end{abstract}
	
\section{Introduction}

Let $\Gamma=(V,E)$ be a strongly regular graph (SRG) with parameters $(n,k,\lambda,\mu)$. For each edge $e=\{u,v\}\in E$, we add one additional vertex $z_e \notin V$ such that it forms 3-uniform hypergraph $\mathcal{H}(\mathcal{V},\mathcal{E})$ where $\mathcal{V}=V \cup Z$, $Z=\{z_e: e \in E\}$ and $\mathcal{E}=\{h_e=\{u,v,z_e\}: u,v \in V, z_e \in Z\}$ \cite{maulana, kostochka}.  Strongly regular graphs are closely related to association schemes. Some examples include the Petersen graph $SRG(10,3,0,1)$ and Paley graph.

Association scheme captures regularity in pairwise relations. It is known that $\operatorname{SRG}(n,k,\lambda,\mu)$ is $2-class$ association scheme which consist of two non trivial relations: adjacency and non adjacency \cite{cameron}. Besides, there are two important examples of association scheme: first, any regular representation of a finite group and second, the distance matrices of any distance-regular graph \cite{hirasaka}. In some cases, association schemes may be too restrictive for complex structures such as hypergraphs with distinct vertex types.

A coherent configuration generalizes an association scheme by allowing multiple fibers and more general relations. In other word, association scheme is a particular case of a coherent configuration with one fiber or we can call homogenoues coherent configuration \cite{klin}. For a permutation group \(G\leq\operatorname{Sym}(X)\), its orbitals on \(X\times X\) form a Schurian coherent configuration, providing a connection between group actions and relational structures \cite{higman}. 

In this paper, we investigate coherent configurations arising from 3-uniform hypergraphs \(\mathcal H(\Gamma)\) constructed from strongly regular graphs \(\Gamma\). Our study focuses on the relationship between their automorphism groups and the induced orbital structures. The Petersen graph is considered as a specific example to examine the basic relations, valencies, intersection numbers, and associated adjacency algebra.

\section{Preliminaries}

\begin{definition}
$\Gamma=(V,E)$ is strongly regular graph with parameters $(n,k,\lambda,\mu)$ where $n=|V|$, every vertex has degree $k$, every adjacent pair has $\lambda$ common neighbors, and every distinct nonadjacent pair has $\mu$ common neighbors.\\ 
We assume $0<k<n-1$. If $A$ is its adjacency matrix, then
\begin{center} 
$A^2=kI+\lambda A+\mu(J-I-A),\qquad (n-k-1)\mu=k(k-\lambda-1)$
\end{center}
In particular $|E|=\frac{nk}{2}$.
\end{definition}

\begin{definition}
The hypergraph associated with $\Gamma$ is $\mathcal{H}(\Gamma)= \mathcal H (\mathcal V,\mathcal E)$, where
\begin{center}
$\mathcal{V}=V \cup Z, Z=\{z_e: e \in E\}$ and $\mathcal{E}= \{h_e=\{u,v,z_e\}: u,v \in V, z_e \in Z\}$.
\end{center}
\end{definition}	
	
\begin{proposition}
For $\Gamma$ with parameters $(n,k,\lambda,\mu)$, the hypergraph $\mathcal H(\Gamma)$ has
\begin{center}
$|\mathcal{V}|=n+nk/2,\quad |\mathcal E|=nk/2,\quad \deg_{\mathcal H}(x)= \begin{cases} k,&x\in V,\\1,&x\in Z\end{cases}$
\end{center}
and it is linear, distinct hyperedges intersect in at most one vertex.
\end{proposition}
\begin{proof}
Each edge \(e=\{u,v\}\) $\in E$ corresponds to a unique hyperedge $h_e=\{u,v,z_e\} \in \mathcal E$. Since \(\Gamma\) is \(k\)-regular, each original vertex in $V$ belongs to exactly \(k\) hyperedges, while each new vertex in $Z$ belongs to exactly one. Moreover, distinct hyperedges share at most one vertex, since their corresponding edges in \(\Gamma\) have at most one common endpoint. Hence, \(\mathcal H(\Gamma)\) is linear.
\end{proof}

\begin{definition}
Let $\mathcal H=(\mathcal V,\mathcal E)$ be a hypergraph. For any two distinct vertices $u,v\in \mathcal V$, the codegree $d_{\mathcal H}$ of $u$ and $v$ is defined as the number of hyperedges containing both $u$ and $v$. 
\begin{center}
$d_{\mathcal H}(u,v)
=
\left|\left\{h=e\in\mathcal E:\{u,v\}\subseteq h_e\right\}\right|.$
\end{center}
\end{definition}

\begin{definition}
An automorphism of $\Gamma$ is a permutation $g\in\Sym(V)$ preserving adjacency. An automorphism of $\mathcal H$ is a permutation of $\mathcal V$ preserving its hyperedge set. If a group $G$ acts on $\mathcal V$, its induced action on ordered pairs is $g\cdot(x,y)=(g x,g y)$. The orbital of $(x,y)$ is $\Orb_G(x,y)=\{(g x,g y):g\in G\}.$
\end{definition}

\begin{definition}
A coherent configuration $(X,\mathcal R)$ is a partition $\mathcal R$ of $X^2$ into nonempty binary relations satisfying:
\begin{enumerate}[label=(CC\arabic*)]
\item $\Delta_X$ is a union of relations in $\mathcal R$;
\item $R^T\in\mathcal R$ for every $R\in\mathcal R$;
\item for $R_i,R_j,R_k\in\mathcal R$, the number
\begin{center}
$p_{ij}^k=|\{z\in X:(x,z)\in R_i,\ (z,y)\in R_j\}|$
\end{center}
is independent of $(x,y)\in R_k$.
\end{enumerate}
A fiber is a set $F\subseteq X$ for which $\Delta_F$ is a basic relation. The configuration is homogeneous if it has one fiber and nonhomogeneous otherwise. It is Schurian if its relations are the orbitals of a permutation group.
\end{definition}

\begin{definition}
An association scheme is a homogeneous coherent configuration. It is symmetric when all its relations are self-transpose, and commutative when its adjacency matrices commute.
\end{definition}

\begin{proposition}
Let $\Gamma=(V,E)$ be a nontrivial strongly regular graph, then the relations defined by
\begin{center}
$R_0=\Delta_V,\qquad R_1=\{(u,v)\in V^2:u\sim v\}, \qquad R_2=\{(u,v)\in V^2:u\neq v,\ u\not\sim v\}$
\end{center}
form a symmetric 2-class association scheme on $V$.
\end{proposition}
\begin{proof}
Let $A$ be the adjacency matrix of $\Gamma$. Define $A_0=I$, $A_1=A$, and $A_2=J-I-A$. \\
These matrices are symmetric and satisfy $A_0+A_1+A_2=J$. \\
By the strongly regular property,
$A_1^2=kA_0+\lambda A_1+\mu A_2$.\\
Furthermore,
\begin{center}
$A_1A_2=(k-\lambda-1)A_1+(k-\mu)A_2$
and\\
$A_2^2=(v-k-1)A_0+(v-2k+\lambda)A_1+(v-2k+\mu-2)A_2$.
\end{center}
Thus, the span of $A_0,A_1,A_2$ is closed under matrix multiplication. The corresponding relations satisfy the association scheme axioms and are symmetric. Hence, $\Gamma$ induces a symmetric 2-class association scheme.
\end{proof}

\section{Results}
\subsection{Automorphisms and the orbital configuration}
\begin{theorem}
For every nontrivial strongly regular graph $\Gamma$,
\[\Aut(\mathcal H(\Gamma))\cong\Aut(\Gamma).\]
\end{theorem}
\begin{proof}
Every $g\in\Aut(\Gamma)$ extends to $\widehat g\in\Sym(\mathcal V)$ by
\begin{center}
$\widehat g(x)=\begin{cases}g(x),&x\in V,\\z_{\{g(u),g(v)\}},&x=z_{\{u,v\}}\in Z.\end{cases}$
\end{center}
This defines an injective homomorphism \(\Aut(\Gamma)\to\Aut(\mathcal H(\Gamma))\). Conversely, since \(k\ge2\), every hypergraph automorphism preserves \(V\) and \(Z\) by their distinct vertex degrees. Its restriction to \(V\) preserves adjacency and thus belongs to \(\Aut(\Gamma)\). Since each new vertex corresponds uniquely to an edge of \(\Gamma\), the automorphism is determined by this restriction. Hence, the homomorphism is surjective and therefore an isomorphism. \(\square\)
\end{proof}

\begin{theorem}
Let $G=\Aut(\mathcal H(\Gamma))$ and let $\mathcal R$ be its orbitals on $\mathcal V \times \mathcal V$. Then $(\mathcal V,\mathcal R)$ is a coherent configuration.
\end{theorem}
\begin{proof}
The orbitals partition $\mathcal V^2$,\\
\emph{(CC1).} For any $g\in G$, $g\cdot(v,v)=(gv,gv)$. Thus, every orbital containing a diagonal pair lies entirely in  $\Delta_\mathcal V$, so $\Delta_\mathcal V$ is a union of orbitals.\\
\emph{(CC2).} If $R=\Orb_G(v,w)$, then $R^T=\Orb_G(w,v)$, which is also an orbital.\\
\emph{(CC3).} Let $R_i,R_j,R_k\in\mathcal R$ and $(v,w),(v',w')\in R_k$. Since both pairs belong to the same orbital, there exists $g\in G$ such that $g(v)=v'$ and $g(w)=w'$. The map $x\mapsto g(x)$ induces a bijection
${x\in \mathcal V:(v,x)\in R_i,\ (x,w)\in R_j}
\longrightarrow
{x\in \mathcal V:(v',x)\in R_i,\ (x,w')\in R_j}.$
Hence, the intersection numbers $p_{ij}^{k}$ are independent of the choice of $(v,w)\in R_k$.\\
Therefore, all three axioms hold, and $(\mathcal V,\mathcal R)$ is a coherent configuration. 
\end{proof}

\begin{corollary}
The coherent configuration $(\mathcal V,\mathcal R)$  is Schurian.
\end{corollary}
\begin{proof}
Its basic relations are precisely the orbitals of $G$.
\end{proof}

\begin{corollary}
The coherent configuration $(\mathcal V,\mathcal R)$ is nonhomogeneous.
\end{corollary}
\begin{proof}
Since vertices in \(V\) and \(Z\) have distinct hypergraph degrees \(k\) and \(1\), respectively, no orbit of \(G\) contains vertices from both sets. Thus, the diagonal contains at least two distinct orbitals, so the coherent configuration is nonhomogeneous.
\end{proof}

\subsection{Case: The Petersen Graph}
\begin{proposition}
Let $P=K(5,2)$ be the Petersen graph, which is an \(\operatorname{SRG}(10,3,0,1)\) with \(\operatorname{Aut}(P)\cong S_5\). The associated hypergraph \(\mathcal H(P)\) has \(25\) vertices and \(15\) hyperedges, and its orbital coherent configuration has two fibers of sizes \(10\) and \(15\), with vertex degrees \(3\) and \(1\), respectively.
\end{proposition}
\begin{proof}
The vertices of $P$ are the ten 2-subsets of $\{1,2,3,4,5\}$, with adjacency defined by disjointness. This gives the parameters $(10,3,0,1)$, and the natural action of $S_5$ yields $\operatorname{Aut}(P)\cong S_5$. Since $P$ has $15$ edges, the construction produces a hypergraph with $25$ vertices and $15$ hyperedges. Moreover, $S_5$ acts transitively on both vertices and edges of $P$, so the associated coherent configuration has exactly two fibers of sizes $10$ and $15$. Their corresponding hypergraph vertex degrees are $3$ and $1$, respectively.
\end{proof}

\begin{proposition}
Let $e,f\in E(P)$ be disjoint edges of the Petersen graph, define
$c(e,f)=\bigl|\{(u,w)\in e\times f:u\sim w\}\bigr|,$
then
$c(e,f)\in\{0,1\}.$
\end{proposition}
\begin{proof}
Let $e=\{u_1,u_2\}$ and $f=\{w_1,w_2\}$ be disjoint edges of $P$. If two cross-edges existed between $e$ and $f$, they would form either a triangle or a 4-cycle together with $e$ and $f$. Both contradict the girth of the Petersen graph, which is 5. Hence, $c(e,f)\leq 1$.
\end{proof}

\begin{theorem}
Under the natural $S_5$-action on $\mathcal V=V\cup Z$, the orbital coherent configuration has exactly $13$ basis relations. They consist of
\begin{enumerate}[label=\textup{(\roman*)}]
 \item the two diagonal relations on $V$ and $Z$;
 \item adjacency and nonadjacency on $V\times V$;
 \item three relations on $V\times Z$, according as $x$ is an endpoint of $e$, or $x\notin e$ with $c(x,e)=0$, or $x\notin e$ with $c(x,e)=1$;
 \item the three transpose relations on $Z\times V$;
 \item three off-diagonal relations on $Z\times Z$: distinct edges sharing an endpoint, disjoint edges with one cross-edge, and disjoint edges with no cross-edge.
\end{enumerate}
\end{theorem}
\begin{proof}
The natural action of $S_5$ on the vertices and edges of the Petersen graph induces an action on $V\cup Z$.\\
On $V\times V$, ordered pairs are classified by their intersection sizes $2, 0,$ and $1$. These correspond to the diagonal, adjacency, and nonadjacency relations, respectively.\\
For $V\times Z$, fix an edge $e=\{A,B\}$. Its stabilizer in $S_5$ partitions $V$ into three orbits: the endpoints of $e$, the vertices adjacent to neither endpoint, and the vertices adjacent to exactly one endpoint. Thus, $V\times Z$ has three orbitals, while $Z\times V$ has their transposes.\\
Finally, on $Z\times Z$, the orbitals are determined by whether two edges are identical, share an endpoint, or are disjoint. By Proposition 4, disjoint edges have either zero or one cross-edge. These cases yield four orbitals.\\
Therefore, the total number of basic relations is\\
$3+3+3+4=13.$
Hence, the orbital coherent configuration has exactly $13$ basic relations.
\end{proof}

\begin{corollary}
Under the indexing
\[
 R_0,R_1,R_2,R_3,R_4,R_5,R_6,R_4^{\mathsf T},R_5^{\mathsf T},R_6^{\mathsf T},R_8,R_9,R_{10},
\]
the cardinalities of the thirteen basic relations are
\[
 10,15,30,60,30,60,60,30,60,60,60,120,30,
\]
respectively. In particular,
\[
 \sum_i |R_i|=625=25^2.
\]
\end{corollary}

\begin{proof}
Since \(P\) has \(10\) vertices and \(15\) edges, the cardinalities on \(V\times V\) follow from its parameters \((10,3,0,1)\). For \(V\times Z\), each edge has two endpoints, while the remaining vertices are classified by their adjacency to the endpoints. The relations on \(Z\times V\) have the same cardinalities as their transposes.

On \(Z\times Z\), each edge meets four other edges and is disjoint from ten edges. Of these ten, eight have exactly one cross-edge and two have none. Thus, the corresponding cardinalities are \(15\), \(60\), \(120\), and \(30\). Summing all thirteen cardinalities gives
\[
\sum_{R\in\mathcal R}|R|=625=|\mathcal V|^2.
\]
\end{proof}

\begin{proposition}
The valencies of the thirteen basic relations are given by the following tables.\[
\begin{array}{c|cccccc}
V & R_0&R_2&R_3&R_4&R_5&R_6\\ \hline
 &1&3&6&3&6&6
\end{array}
\]
and
\[
\begin{array}{c|ccccccc}
Z & R_1&R_4^{\mathsf T}&R_5^{\mathsf T}&R_6^{\mathsf T}&R_8&R_9&R_{10}\\ \hline
 &1&2&4&4&4&8&2.
\end{array}
\]
In particular, the sum of each row is $25$.
\end{proposition}

\begin{proof}
For each basic relation $R\subseteq X\times Y$, where $X,Y\in{V,Z}$, its valency is $n_R=|R|/|X|$. Using the cardinalities from Corollary 3, with $|V|=10$ and $|Z|=15$, yields the stated valencies. In particular, their sums are
\begin{center}
$[
1+3+6+3+6+6=25],
\qquad
[1+2+4+4+4+8+2=25
]$.
\end{center}
\end{proof}

\begin{theorem}
The thirteen relations above form a nonhomogeneous Schurian coherent configuration with exactly two fibers $V$ and $Z$.
\end{theorem}
\begin{proof}
By Theorem 3, the thirteen basic relations are precisely the orbitals of \(G=\operatorname{Aut}(\mathcal H(P))\cong S_5\) on \(\mathcal V\times\mathcal V\). Hence, they form a Schurian coherent configuration. Since \(S_5\) acts transitively on \(V\) and \(Z\) separately, and these sets are distinct orbits, the configuration has exactly two fibers. Therefore, it is nonhomogeneous. \
\end{proof}

\begin{proposition}
For the Petersen orbital coherent configuration every intersection number $p_{ij}^k$ is well-defined. With the index convention
\[
R_0\leftrightarrow0,\ R_1\leftrightarrow1,\ R_2\leftrightarrow2,\ R_3\leftrightarrow3,
\ R_4\leftrightarrow4,\ R_5\leftrightarrow5,\ R_6\leftrightarrow6,
\]
\[
R_4^{\mathsf T}\leftrightarrow7,\quad
R_5^{\mathsf T}\leftrightarrow8,\quad
R_6^{\mathsf T}\leftrightarrow9,\quad
R_8\leftrightarrow10,\quad R_9\leftrightarrow11,\quad R_{10}\leftrightarrow12,
\]
there are exactly $159$ nonzero intersection numbers among the $13^3=2197$ possible triples $(i,j,k)$.  
\end{proposition}

\begin{proof}
Fix \(R_k\) and a representative pair \((x,y)\in R_k\). Since the relations are \(S_5\)-orbitals, the number
\[
p_{ij}^{k}
=
\left|\{z\in\mathcal V:(x,z)\in R_i,\ (z,y)\in R_j\}\right|
\]
is independent of the chosen representative. Counting the intermediate vertices \(z\) for each relation yields exactly \(159\) nonzero intersection numbers, as listed in Table 1. \end{proof}
\begin{longtable}{@{}c c r @{\hspace{1.2em}} c c r@{}}
\caption{All nonzero intersection numbers of the Petersen configuration.}\label{tab:all-pijk}\\
\toprule
Relation & Intersection number & Value & Relation & Intersection number & Value\\
\midrule
\endfirsthead
\toprule
Relation & Intersection number & Value & Relation & Intersection number & Value\\
\midrule
\endhead
$R_{0}$ & $p^{0}_{0,0}$ & $1$ & $R_{7}$ & $p^{7}_{8,3}$ & $4$ \\
$R_{0}$ & $p^{0}_{2,2}$ & $3$ & $R_{7}$ & $p^{7}_{9,2}$ & $2$ \\
$R_{0}$ & $p^{0}_{3,3}$ & $6$ & $R_{7}$ & $p^{7}_{9,3}$ & $2$ \\
$R_{0}$ & $p^{0}_{4,7}$ & $3$ & $R_{7}$ & $p^{7}_{10,7}$ & $2$ \\
$R_{0}$ & $p^{0}_{5,8}$ & $6$ & $R_{7}$ & $p^{7}_{10,9}$ & $2$ \\
$R_{0}$ & $p^{0}_{6,9}$ & $6$ & $R_{7}$ & $p^{7}_{11,8}$ & $4$ \\
$R_{1}$ & $p^{1}_{1,1}$ & $1$ & $R_{7}$ & $p^{7}_{11,9}$ & $4$ \\
$R_{1}$ & $p^{1}_{7,4}$ & $2$ & $R_{7}$ & $p^{7}_{12,8}$ & $2$ \\
$R_{1}$ & $p^{1}_{8,5}$ & $4$ & $R_{8}$ & $p^{8}_{1,8}$ & $1$ \\
$R_{1}$ & $p^{1}_{9,6}$ & $4$ & $R_{8}$ & $p^{8}_{7,3}$ & $2$ \\
$R_{1}$ & $p^{1}_{10,10}$ & $4$ & $R_{8}$ & $p^{8}_{8,0}$ & $1$ \\
$R_{1}$ & $p^{1}_{11,11}$ & $8$ & $R_{8}$ & $p^{8}_{8,2}$ & $1$ \\
$R_{1}$ & $p^{1}_{12,12}$ & $2$ & $R_{8}$ & $p^{8}_{8,3}$ & $2$ \\
$R_{2}$ & $p^{2}_{0,2}$ & $1$ & $R_{8}$ & $p^{8}_{9,2}$ & $2$ \\
$R_{2}$ & $p^{2}_{2,0}$ & $1$ & $R_{8}$ & $p^{8}_{9,3}$ & $2$ \\
$R_{2}$ & $p^{2}_{2,3}$ & $2$ & $R_{8}$ & $p^{8}_{10,8}$ & $2$ \\
$R_{2}$ & $p^{2}_{3,2}$ & $2$ & $R_{8}$ & $p^{8}_{10,9}$ & $2$ \\
$R_{2}$ & $p^{2}_{3,3}$ & $4$ & $R_{8}$ & $p^{8}_{11,7}$ & $2$ \\
$R_{2}$ & $p^{2}_{4,7}$ & $1$ & $R_{8}$ & $p^{8}_{11,8}$ & $2$ \\
$R_{2}$ & $p^{2}_{4,9}$ & $2$ & $R_{8}$ & $p^{8}_{11,9}$ & $4$ \\
$R_{2}$ & $p^{2}_{5,8}$ & $2$ & $R_{8}$ & $p^{8}_{12,7}$ & $1$ \\
$R_{2}$ & $p^{2}_{5,9}$ & $4$ & $R_{8}$ & $p^{8}_{12,8}$ & $1$ \\
$R_{2}$ & $p^{2}_{6,7}$ & $2$ & $R_{9}$ & $p^{9}_{1,9}$ & $1$ \\
$R_{2}$ & $p^{2}_{6,8}$ & $4$ & $R_{9}$ & $p^{9}_{7,2}$ & $1$ \\
$R_{3}$ & $p^{3}_{0,3}$ & $1$ & $R_{9}$ & $p^{9}_{7,3}$ & $1$ \\
$R_{3}$ & $p^{3}_{2,2}$ & $1$ & $R_{9}$ & $p^{9}_{8,2}$ & $2$ \\
$R_{3}$ & $p^{3}_{2,3}$ & $2$ & $R_{9}$ & $p^{9}_{8,3}$ & $2$ \\
$R_{3}$ & $p^{3}_{3,0}$ & $1$ & $R_{9}$ & $p^{9}_{9,0}$ & $1$ \\
$R_{3}$ & $p^{3}_{3,2}$ & $2$ & $R_{9}$ & $p^{9}_{9,3}$ & $3$ \\
$R_{3}$ & $p^{3}_{3,3}$ & $3$ & $R_{9}$ & $p^{9}_{10,7}$ & $1$ \\
$R_{3}$ & $p^{3}_{4,8}$ & $2$ & $R_{9}$ & $p^{9}_{10,8}$ & $2$ \\
$R_{3}$ & $p^{3}_{4,9}$ & $1$ & $R_{9}$ & $p^{9}_{10,9}$ & $1$ \\
$R_{3}$ & $p^{3}_{5,7}$ & $2$ & $R_{9}$ & $p^{9}_{11,7}$ & $2$ \\
$R_{3}$ & $p^{3}_{5,8}$ & $2$ & $R_{9}$ & $p^{9}_{11,8}$ & $4$ \\
$R_{3}$ & $p^{3}_{5,9}$ & $2$ & $R_{9}$ & $p^{9}_{11,9}$ & $2$ \\
$R_{3}$ & $p^{3}_{6,7}$ & $1$ & $R_{9}$ & $p^{9}_{12,9}$ & $2$ \\
$R_{3}$ & $p^{3}_{6,8}$ & $2$ & $R_{10}$ & $p^{10}_{1,10}$ & $1$ \\
$R_{3}$ & $p^{3}_{6,9}$ & $3$ & $R_{10}$ & $p^{10}_{7,4}$ & $1$ \\
$R_{4}$ & $p^{4}_{0,4}$ & $1$ & $R_{10}$ & $p^{10}_{7,6}$ & $1$ \\
$R_{4}$ & $p^{4}_{2,4}$ & $1$ & $R_{10}$ & $p^{10}_{8,5}$ & $2$ \\
$R_{4}$ & $p^{4}_{2,6}$ & $2$ & $R_{10}$ & $p^{10}_{8,6}$ & $2$ \\
$R_{4}$ & $p^{4}_{3,5}$ & $4$ & $R_{10}$ & $p^{10}_{9,4}$ & $1$ \\
$R_{4}$ & $p^{4}_{3,6}$ & $2$ & $R_{10}$ & $p^{10}_{9,5}$ & $2$ \\
$R_{4}$ & $p^{4}_{4,1}$ & $1$ & $R_{10}$ & $p^{10}_{9,6}$ & $1$ \\
$R_{4}$ & $p^{4}_{4,10}$ & $2$ & $R_{10}$ & $p^{10}_{10,1}$ & $1$ \\
$R_{4}$ & $p^{4}_{5,11}$ & $4$ & $R_{10}$ & $p^{10}_{10,10}$ & $1$ \\
$R_{4}$ & $p^{4}_{5,12}$ & $2$ & $R_{10}$ & $p^{10}_{10,11}$ & $2$ \\
$R_{4}$ & $p^{4}_{6,10}$ & $2$ & $R_{10}$ & $p^{10}_{11,10}$ & $2$ \\
$R_{4}$ & $p^{4}_{6,11}$ & $4$ & $R_{10}$ & $p^{10}_{11,11}$ & $4$ \\
$R_{5}$ & $p^{5}_{0,5}$ & $1$ & $R_{10}$ & $p^{10}_{11,12}$ & $2$ \\
$R_{5}$ & $p^{5}_{2,5}$ & $1$ & $R_{10}$ & $p^{10}_{12,11}$ & $2$ \\
$R_{5}$ & $p^{5}_{2,6}$ & $2$ & $R_{11}$ & $p^{11}_{1,11}$ & $1$ \\
$R_{5}$ & $p^{5}_{3,4}$ & $2$ & $R_{11}$ & $p^{11}_{7,5}$ & $1$ \\
$R_{5}$ & $p^{5}_{3,5}$ & $2$ & $R_{11}$ & $p^{11}_{7,6}$ & $1$ \\
$R_{5}$ & $p^{5}_{3,6}$ & $2$ & $R_{11}$ & $p^{11}_{8,4}$ & $1$ \\
$R_{5}$ & $p^{5}_{4,11}$ & $2$ & $R_{11}$ & $p^{11}_{8,5}$ & $1$ \\
$R_{5}$ & $p^{5}_{4,12}$ & $1$ & $R_{11}$ & $p^{11}_{8,6}$ & $2$ \\
$R_{5}$ & $p^{5}_{5,1}$ & $1$ & $R_{11}$ & $p^{11}_{9,4}$ & $1$ \\
$R_{5}$ & $p^{5}_{5,10}$ & $2$ & $R_{11}$ & $p^{11}_{9,5}$ & $2$ \\
$R_{5}$ & $p^{5}_{5,11}$ & $2$ & $R_{11}$ & $p^{11}_{9,6}$ & $1$ \\
$R_{5}$ & $p^{5}_{5,12}$ & $1$ & $R_{11}$ & $p^{11}_{10,10}$ & $1$ \\
$R_{5}$ & $p^{5}_{6,10}$ & $2$ & $R_{11}$ & $p^{11}_{10,11}$ & $2$ \\
$R_{5}$ & $p^{5}_{6,11}$ & $4$ & $R_{11}$ & $p^{11}_{10,12}$ & $1$ \\
$R_{6}$ & $p^{6}_{0,6}$ & $1$ & $R_{11}$ & $p^{11}_{11,1}$ & $1$ \\
$R_{6}$ & $p^{6}_{2,4}$ & $1$ & $R_{11}$ & $p^{11}_{11,10}$ & $2$ \\
$R_{6}$ & $p^{6}_{2,5}$ & $2$ & $R_{11}$ & $p^{11}_{11,11}$ & $4$ \\
$R_{6}$ & $p^{6}_{3,4}$ & $1$ & $R_{11}$ & $p^{11}_{11,12}$ & $1$ \\
$R_{6}$ & $p^{6}_{3,5}$ & $2$ & $R_{11}$ & $p^{11}_{12,10}$ & $1$ \\
$R_{6}$ & $p^{6}_{3,6}$ & $3$ & $R_{11}$ & $p^{11}_{12,11}$ & $1$ \\
$R_{6}$ & $p^{6}_{4,10}$ & $1$ & $R_{12}$ & $p^{12}_{1,12}$ & $1$ \\
$R_{6}$ & $p^{6}_{4,11}$ & $2$ & $R_{12}$ & $p^{12}_{7,5}$ & $2$ \\
$R_{6}$ & $p^{6}_{5,10}$ & $2$ & $R_{12}$ & $p^{12}_{8,4}$ & $2$ \\
$R_{6}$ & $p^{6}_{5,11}$ & $4$ & $R_{12}$ & $p^{12}_{8,5}$ & $2$ \\
$R_{6}$ & $p^{6}_{6,1}$ & $1$ & $R_{12}$ & $p^{12}_{9,6}$ & $4$ \\
$R_{6}$ & $p^{6}_{6,10}$ & $1$ & $R_{12}$ & $p^{12}_{10,11}$ & $4$ \\
$R_{6}$ & $p^{6}_{6,11}$ & $2$ & $R_{12}$ & $p^{12}_{11,10}$ & $4$ \\
$R_{6}$ & $p^{6}_{6,12}$ & $2$ & $R_{12}$ & $p^{12}_{11,11}$ & $4$ \\
$R_{7}$ & $p^{7}_{1,7}$ & $1$ & $R_{12}$ & $p^{12}_{12,1}$ & $1$ \\
$R_{7}$ & $p^{7}_{7,0}$ & $1$ & $R_{12}$ & $p^{12}_{12,12}$ & $1$ \\
$R_{7}$ & $p^{7}_{7,2}$ & $1$ &  & &  \\
\bottomrule
\end{longtable}

\begin{proposition}
Let \(A_0,\ldots,A_{12}\) be the adjacency matrices of the thirteen basic relations of the Petersen orbital coherent configuration. Then, for all \(0\leq i,j\leq12\),
\[
A_iA_j=\sum_{k=0}^{12}p_{ij}^{k}A_k.
\]
\end{proposition}
\begin{proof}
For any \(x,y\in\mathcal V\), the \((x,y)\)-entry of \(A_iA_j\) counts the vertices \(z\in\mathcal V\) satisfying \((x,z)\in R_i\) and \((z,y)\in R_j\). If \((x,y)\in R_k\), this number equals \(p_{ij}^{k}\). Since the basic relations partition \(\mathcal V\times\mathcal V\), the identity follows by comparing corresponding matrix entries.
\end{proof}

\begin{acknowledgement}
The author is grateful to Professor Manabu Oura (Kanazawa University) for his valuable suggestions and helpful discussions.
\end{acknowledgement}

\end{document}